\documentclass[11pt]{article}
\usepackage[
	a4paper,
	margin=1in
]{geometry}
\usepackage{amsmath, amssymb, amsfonts, amsthm, graphicx, fullpage, subcaption,enumitem}
\usepackage{hyperref}
\usepackage{cleveref}

\usepackage[utf8]{inputenc} 
\usepackage[T1]{fontenc}

\usepackage{hyperref}

\hypersetup{
    colorlinks=true,
    linkcolor=blue,
    filecolor=magenta,      
    urlcolor=cyan,
}

\theoremstyle{plain} \newtheorem{theorem}{Theorem}
\theoremstyle{plain} \newtheorem{problem}[theorem]{Problem}
\theoremstyle{plain} \newtheorem{lemma}{Lemma}[section]
\theoremstyle{plain} 
\theoremstyle{plain} \newtheorem{claim}{Claim}[section]
\theoremstyle{definition} \newtheorem{definition}{Definition}[section]
\theoremstyle{plain}

\title{On the size of universal graphs for spanning trees}
\author{Jaehoon Kim\thanks{Department of Mathematical Sciences, KAIST, South Korea.
        \emph{E-mails:} \textbf{\{jaehoon.kim, minseo.kim00\}@kaist.ac.kr}}~\thanks{Supported by the National Research Foundation of Korea (NRF) grant funded by the Korea government(MSIT) No. RS-2023-00210430.}
\and  Minseo Kim\footnotemark[1]}
\date{}

\begin{document}

\maketitle

\begin{abstract}
Chung and Graham~[J. London Math. Soc., 1983]  claimed that there exists an $n$-vertex graph $G$ containing all $n$-vertex trees as subgraphs that has at most \( \frac{5}{2}n \log_2 n + O(n)\) edges. 
We identify an error in their proof. This error can be corrected by adding more edges, which increases the number of edges to \( e(G) \le \frac{7}{2}n \log_2 n + O(n). \) Moreover, we further improve this by showing that there exists such an $n$-vertex graph with at most \(
\left(5- \frac{1}{3}\right)n \log_3 n + O(n) \leq  2.945 \cdot n\log_2(n)\) edges. This is the first improvement of the bound since Chung and Graham’s pioneering work four decades ago.
\end{abstract}

\section{introduction}

For a graph family $\mathcal{F}$, an $\mathcal{F}$-universal graph is a graph which contains all elements of $\mathcal{F}$ as subgraphs. 
The minimum number of edges in $\mathcal{F}$-universal graphs has been extensively studied for various graph families $\mathcal{F}$.
To list a few, 
the family of all bounded degree graphs \cite{alonsecond, szemeredi, alonthird, capalbosecond}, the family of all $m$-edge graphs \cite{alon, babai, bucic}, the families of planar graphs with certain restrictions \cite{bhatt, capalbo, joret}, the families of separable bounded degree graphs \cite{chung} and the family of trees \cite{grahamfirst, grahamsecond, grahamthird, pipenger, gyori, montgomery}.

Chung and Graham \cite{grahamfirst, grahamthird} considered the natural problem of determining the minimum number of edges of a $\mathcal{T}_n$-universal graph, where $\mathcal{T}_n$ is the family of all $n$-vertex trees. They define $s^*(n)$ as the minimum number of edges in such a graph. Additionally, they define $s(n)$ as the minimum number of edges in $n$-vertex $\mathcal{T}_n$-universal graphs. In other words, $s(n)$ is the minimum number of edges in an $n$-vertex graph that contains all possible spanning trees. 

Chung and Graham \cite{grahamfirst} proved a lower bound on $s^*(n)$ and Chung and Graham \cite{grahamthird} claimed an upper bound on $s(n)$ yielding the following bound.
\begin{equation}\label{eqn}
\frac{1}{2}n\ln{n}< s^{*}(n)\leq s(n)\leq \frac{5}{2}n\log_2{n}+O(n).
\end{equation}
Unfortunately, we found an error in the proof of the upper bound in \eqref{eqn}. However, this error can be fixed by adding more edges to their construction and it yields the following weaker upper bound.
\begin{theorem}\label{thm1} \(s(n)\leq\frac{7}{2}n\log_2{n}+ 4n.\)
\end{theorem}
We explain this mistake and how to correct it in \Cref{ch3}.
We note that this mistake was independently discovered by Kaul and Wood in \cite{kaul}, who provided an alternative weaker bound $s(n)= O(n\log n \log\log n)$. In fact, Bhatt, Chung, Leighton, and Rosenberg \cite{bhatt} and Chung \cite{chung} both cited the upper bound from \cite{grahamthird} as $\frac{7}{2}n\log_2 n + O(n)$. So we believe that they knew about the mistake and the way to fix it as we present in \Cref{ch3}.

Recently, Gy\H ori, Li, Salia and Tompkins in \cite{gyori} improved the lower bound to \(n\ln{n}-O(n)\leq s^{*}(n).\)
So the state of the art is \[n\ln{n}-O(n)\leq s^{*}(n)\leq s(n)\leq \frac{7}{2}n\log_2{n}+O(n).\]
Here, the gap between two bounds are a multiplication factor of a number $\frac{7}{2\ln(2)} = 5.04\dots$, slightly larger than $5$.
On the other hand, the upper bound has remained unchanged since \cite{grahamthird} for over four decades. In this paper, we present the first improvement since \cite{grahamthird}.

\begin{theorem} \label{thm2}
\(s(n)\leq \frac{14}{3}n\log_3{n} + O(n).\)
\end{theorem}

Note that $\frac{14}{3\ln(3)}$ is roughly $4.247\dots$, thus the gap between the two best known bounds are slightly below $4.25$.

\section{Preliminaries} \label{ch2}

In this section, we collect some useful definitions and lemmas. Many of the definitions in this section are modifications of the ones from \cite{grahamthird}.

For a given graph $H$ and $G$, we say $H$ {\bf embeds} into $G$ if $G$ contains a subgraph isomorphic to $H$, and $\phi:V(H)\rightarrow V(G)$ is an {\bf embedding} if it is an isomorphism from $H$ to a corresponding subgraph $\phi(H)$ of $G$. We write $|H|$ to denote the number of vertices in $H$ and $N_H(x)$ denotes the neighborhood of $x$ in $H$. We write $H-U$ to denote the graph obtained from $H$ by deleting vertices in $U$. Sometimes we write $G-H$ to denote $G-V(H)$ for a subgraph $H$ of $G$. We say a graph $G$ is {\bf universal} if it contains all $|G|$-vertex trees as subgraphs.

For an $n$-vertex rooted tree $T$ with the root $v$, the {\bf level} $L_T(u)$ of a vertex $u\in V(T)$ is the distance between $u$ and $v$. The level $L(T)$ of a rooted tree $T$ is $\max \{L_T(u): u\in V(T)\}$. Let $D_T(u)$ be the set of all descendants of $u$ and let  $D_T[u]=D_T(u)\cup \{u\}$. Sometimes, we write $D_T[u]$ to denote the subtree of $T$ induced by the vertex set $D_T[u]$ if it is clear from the context. We also write $\nu_T(u)= |D_T[u]|$. We omit the subscript $T$ in these notations if it is clear from the context. 

For each vertex $u$, we give an order on the children from left to right. We perform DFS(Depth First Search) starting from the root $v$, visiting the left-most unvisited child of the current vertex. We backtrack to the parent if all children of the current vertex have been visited. This DFS produces a natural ordering $x_1,\dots, x_n$ of the vertices in $V(T)$ in the ordering of the first appearance of the vertex in this DFS. This order is called {\bf the DFS preorder traversal}. In particular, the root $v$ is $x_1$.
We say that a vertex set $S\subseteq V(T)$ is {\bf $T$-admissible} if $S$ is the set $\{x_1,\dots, x_k\}$ of the initial $k$ vertices in DFS preorder traversal for some $k$. In particular, the empty set and the one vertex set $\{v\}=\{x_1\}$ are also $T$-admissible.

For a vertex $u$, consider the unique path $u_0, u_1,\dots, u_t$ in $T$ from $u=u_0$ to the root $v=u_t$. The vertex $u_i$ is the {\bf $i$-th ancestor} of $u$ and $u$ is an {\bf $i$-th descendant} of $u_i$.  In particular, the $1$-st ancestor is called the {\bf parent} and a $1$-st descendant is called a {\bf child}.  For all $i\geq t$, we let the {\bf $i$-th ancestor} be the root $v$ for our convenience. We write $u^*$ to denote the parent of $u$. {\bf Siblings} of $u$ are vertices sharing the same parent with $u$ and {\bf cousins} of $u$ are  vertices having the same level as $u$. {\bf Left-siblings}/{\bf left-cousins} of $u$ are the siblings/cousins of $u$ that appear before $u$ in the ordering $x_1,\dots, x_n$. Similarly, we define {\bf right-siblings} and {\bf right-cousins} as well.
The {\bf nearest-left cousin} $l(u)$ of $u=x_k$ is $x_i$ with the maximum $i<k$ such that $L_T(x_i)=L_T(x_k)$. In other words, it is the vertex on the same level as $u$ that is to the left of $u$ but closest to $u$ in the order $x_1,\dots, x_n$. In particular, it is a sibling of $u$ if $u$ is not the left-most child of its parent. Otherwise, it is the right-most child of the nearest-left cousin of its parent, assuming it exists.

For a given rooted tree $(T,v)$, we generate a supergraph of $T$ with the same vertex set as follows. 
\begin{definition}\label{newdef}
For a given rooted tree $(T,v)$, {\bf $\overrightarrow{G}^{r}_T$} is a directed graph on the vertex set $V(T)$ that is obtained from $T$ by adding all directed edges as follows for every $u\in V(T).$
\begin{enumerate}[label=\upshape \textbf{G\arabic{enumi}}]
     \item \label{G1} $\overrightarrow{uw}$ for all descendants $w$ of $u$.
     \item \label{G2} $\overrightarrow{uw}$ for all left-siblings $w$ of $u$, and $\overrightarrow{uw'}$ for all descendants $w'$ of $w$.
     \item \label{G3} $\overrightarrow{uw}$ for nearest-left cousin $w=l(u^*)$ of the parent $u^*$ of $u$, and $\overrightarrow{uw'}$ for all descendants $w'$ of $w$.
     \item \label{G4} $\overrightarrow{uw}$ for any $i$-th descendant $w$ of $u'$ for each $i\in [r]$ where $u'$ is either the $r$-th ancestor of $u$ or the nearest-left cousin of the $r$-th ancestor of $u$.
\end{enumerate}
For the directed graph $\overrightarrow{G}^{r}_T$, the underlying undirected simple graph is denoted by {\bf $G^{r}_T$} and it is called  the graph {\bf$r$-generated} by $T$. For $r=0$, we omit the superscript $r$ and say that $G_T$ is generated by $T$.
\end{definition}
Note that \ref{G4} does not add any additional directed edges if $r\in \{0,1\}$. We say that an induced subgraph $G^r_T[S]\subseteq G^r_T$ is {\bf $T$-admissible} if $S$ is a $T$-admissible set. 
An embedding $\phi: T'\rightarrow G^{r}_T$ is {\bf a $T$-admissible embedding} if the set $V(G^{r}_T)\setminus \phi(T')$ is $T$-admissible.
Given a forest $T$ and a vertex $u\in V(T)$, a component in $T-u$ is called a {\bf $u$-component of $T$}. In fact, the following simple lemma about $u$-components is useful for us. We include a proof for completeness. For sets $C_1,\dots, C_t$ and $I\subseteq [t]$, we write $C_I$ to denote $\bigcup_{i\in I} C_i$.

\begin{lemma} \label{lem2}
Let $x\in \mathbb{N}$ and let $T$ be a forest on at least $x+1$ vertices.
For each $u\in V(T)$, there exists a vertex $w$ and a collection $\{C_1,\dots, C_t\}$ of $w$-components satisfying $u\notin C_{[t]}$ and
\[x\leq \left| C_{[t]}\right|\leq 2x-1.\]
\end{lemma}
\begin{proof}
By adding edges if necessary, we may assume that $T$ is a tree.
Consider $u$-components $C_1,\dots, C_p$. If each of them contains less than $2x$ vertices, we can easily choose a collection as desired with $w=u$.  Otherwise, assume $C_1$ has at least $2x$ vertices and $u_1$ is the unique neighbor of $u$ in $C_1$.
Consider $u_1$ components $D_1,\dots, D_{s}$ not containing $u$. Again, if each of them contains less than $2x$ vertices, we can obtain a desired collection with $w=u_1$. Otherwise, assume $D_1$ has at least $2x$ vertices and $u_2$ is the unique  neighbor of $u_1$ in $D_1$. We repeat this for $u_2,u_3,\dots$, then we obtain a desired collection when the process stop. As the graph is a finite tree, it is obvious that the process must stop. This proves the lemma.
\end{proof}

\begin{definition}
    We say that a collection $C_1,\dots, C_t$ of $w$-components in a forest $T$ is {\bf ${\bf (x,y)}$-feasible} if $x\leq |C_{[t]}\cup \{w\} |\leq x+y-2$. We say that the collection is {\bf ${\bf (x,y)}$-critical} if $x+y-2 \leq  |C_{[t]}| \leq  2x-3$ and $|C_{I}|\leq  x-2$ for any $I\subsetneq [t]$ and $t\geq 2$.
\end{definition}
 In particular, we know that $(x,y)$-critical collection satisfies $|C_i|\geq y$ for $i\in [t]$. Also, \Cref{lem2} ensures that $(x,x)$-feasible collection always exists in $T$ for $x\geq 2$ when $T$ has at least $x+1$ vertices. Then the following claim guarantees the existence of either an $(x,y)$-feasible collection or an $(x,y)$-critical collection.
\begin{lemma} \label{lem: feasible critical}
    For a forest $T$ on at least $x+1$ vertices and $x>y\geq 2$. For a vertex $u\in V(T)$, there exists a vertex $w$ and a collection $\{C_1,\dots, C_t\}$ of $w$-components satisfying $u\notin C_{[t]}$ that is either $(x,y)$-feasible or $(x,y)$-critical.
\end{lemma}
\begin{proof}
Suppose that such an $(x,y)$-feasible collection does not exist. Then we can assume that $|T|> x+y-2$, as otherwise the collection of all $u$-components is $(x,y)$-feasible.

   Suppose we have a vertex $w_i$ and a collection of $w_i$-components $C_1,\dots, C_{t}$ not containing $u$ that satisfies $x+y-2\leq |C_{[t]}|\leq 2x-3$.   
Indeed, such a vertex $w_0$ exists, as  \Cref{lem2} implies that there exist a vertex $w_0$ and $w_0$-components $C_1,\dots, C_{t}$ not containing $u$ that satisfies $x-1\leq |C_{[t]}|\leq 2x-3$. As we assume that such an $(x,y)$-feasible collection does not exist, $x+y-2\leq |C_{[t]}|\leq 2x-3$ holds. 
   For given such a vertex $w_i$, we will find $w_{i+1}$ as follows.

    Assume $|C_1|\geq \dots \geq |C_{t}|$ and $s= \max\{ t': |C_{t'}|\geq y \}$. We must have $|C_{[s]}|\geq x+y-2$. If not, we obtain an $(x,y)$-feasible collection $\{C_1,\dots, C_{t'}\}$ for some $t'\geq s$ because $|C_{[t]}|,|C_{[t-1]}|,\dots, |C_{[s]}|$ is a decreasing sequence of difference at most $y-1$ starting at a number at least $x+y-2$ and ending at a number at most $x+y-3$. It must hit a number in the interval $[x-1,x+y-3]$ of $y-1$ consecutive numbers at some point. 
    Take a minimal subcollection $\{D_1,\dots, D_{a_i}\}$ of $\{C_1,\dots, C_s\}$ whose collective size is at least $x+y-2$. If $a_i\geq 2$, then by minimality of $\{D_1,\dots, D_{a_i}\}$, this is a desired $(x,y)$-critical collection with $w_i$. Otherwise, if $a_i=1$ then we let $w_{i+1}$ be the unique neighbor of $w_i$ in $D_1$ and consider $w_{i+1}$-components that lying entirely in $D_{1}$ (thus not containing $u$) to repeat the process. As $T$ is a finite forest, this process must end with a desired $(x,y)$-critical collection. This proves the lemma.
\end{proof}

\subsection{Proof overview}\label{sec: overview}
As explained in \Cref{ch3}, Chung and Graham \cite{grahamthird} considered a graph generated by a complete binary tree $B(k)$ and obtained a universal graph with $\frac{7}{2}n\log_2 n + O(n)$ edges. One natural idea for further improvement is to consider graphs generated by other complete $t$-ary trees for larger values of $t$. It turns out that complete ternary tree yields a slightly better example with $5n\log_3 n + O(n)  \leq  3.156 \cdot  n\log_2n$ edges while all $t>3$ yield worse examples. 

\begin{figure}[h]
\centering
\includegraphics[width=0.5\textwidth]{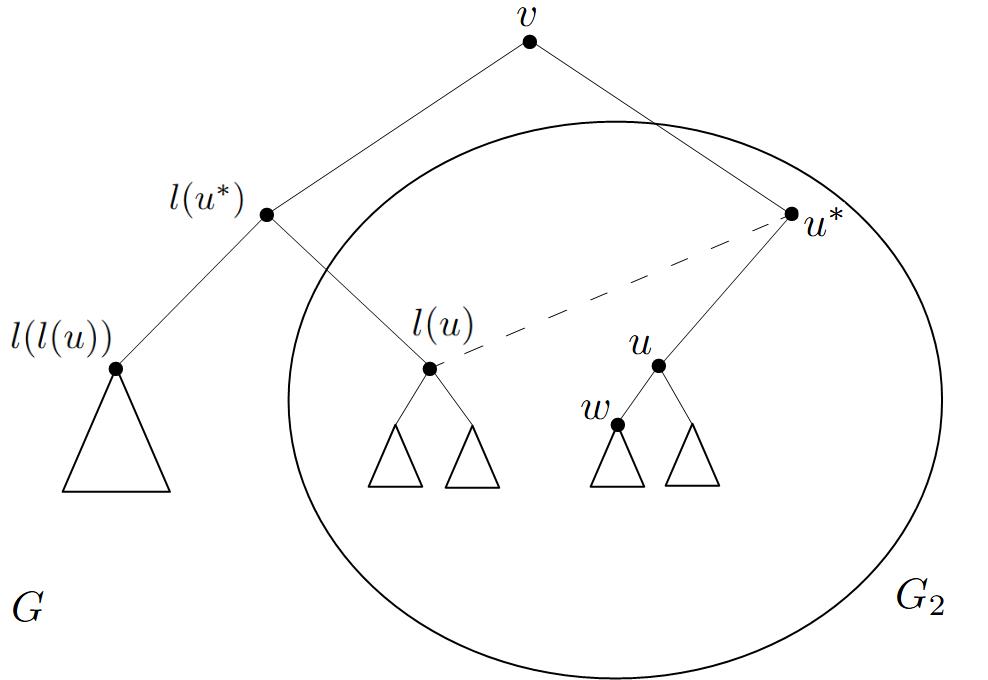}
\caption{ A $G(k)$-admissible graph $G$ } \label{fig1}
\end{figure}

To explain the rough idea in \cite{grahamthird}, suppose that we want to embed a tree $T'$ into $G_T$, where $G_T$ looks like \Cref{fig1}. Then one can use \Cref{lem2} to find a vertex $w$ and some $w$-components $C_1,\dots, C_t$ of $T'$ such that $C_{[t]}$ has size between $\nu(u)$ and $2\nu(u)$. 
Take a subgraph $G_2$ induced by $D_{T}[l(u)]\cup D_{T}[u]\cup \{u^*\}$ as in \Cref{fig1}. By showing this is isomorphic to a $T$-admissible graph, we can use the induction hypothesis to find a $T$-admissible embedding $\phi$ of $T'_0 = T'[\{w\}\cup C_{[t]}]$ into $G_2$ as desired. Furthermore, we ensure that $w$ maps to $u$. Once again, we apply the induction hypothesis to find an embedding $\phi'$ of $T'_1= T'\setminus T'_0$ into $G_{T} - \phi(T'_0)$.
 As all the edges between $T'_0$ and $T'_1$ are incident to $w$ and $u=\phi(w)$ is adjacent to every vertex in $G_T$, we can concatenate $\phi$ and $\phi'$ to obtain a desired embedding. One crucial condition for this argument to work is that $T'_0$ has fewer than $|G_2|$ vertices. This condition is implied by $\nu(l(u))\geq \nu(u)$. Thus, this argument requires that $l(u)$ has at least as many descendants as $u$ for every vertex $u$.

On the other hand, $\nu(l(u))$ being smaller than $\nu(u)$ helps us to obtain a construction with fewer edges. Roughly speaking, this is because the edges we obtain from \Cref{newdef}~\ref{G2} and \ref{G3} become smaller if left siblings (or nearest left-cousin) of $u$ have smaller number of descendants. 
If we allow the ratio between $\nu(l(u))$ and $\nu(u)$ to be as large as $K=\Omega(\varepsilon^{-3})$, then we can construct a tree $T$ with $e(G^r_T)\leq (\frac{5}{2}+\varepsilon) n \log_2 n$. However, we were not able to prove that such a graph contains all possible spanning trees, as the imbalance imposes additional difficulty to embed an arbitrary tree $T'$ into $G^r_T$.

The crucial point is that, in the induction process, $C_{[t]}\cup w$ might be too large to embed into the union of two trees $D_T[l(u)] \cup D_T[u]$. 
So, we sometimes have to merge three trees $D_T[l(l(u))], D_T[l(u)], D_T[u]$ (instead of two) into one to invoke the induction hypothesis.
However, as some edges between $D_T[u]$ and $D_T[l(l(u))]$ might not be present, embedding the tree becomes much more difficult, if not impossible.  
In fact, for large imbalance $K$, we are able to find a desired embedding only by adding additional edges to $G^r_T$. On top of the \Cref{newdef}, if we additionally add edges from each vertex $u$ to its neareset-right cousin and all of its descendants in $G^r_T$, a technical analysis yields that such additions produce a universal graph. Moreover, we can construct a such tree $T$ with  $e(G^r_T)\leq (3+\varepsilon) n \log_2 n$. Although this is a significant improvements from $\frac{7}{2}n\log_2 n+O(n)$ bound by Chung and Graham, we do not present this construction because there is a simpler and better construction.

What we present in this paper is another less technical way to obtain a better upper bound by taking $K=2$. 
The main issue for embedding the subtree $T'_0 = T'[\{w\}\cup C_{[t]}]$ into the union of three trees $D_T[l(l(u))], D_T[l(u)], D_T[u]$ is that 
the edges $ww_1, ww_2,\dots, w w_t$ between the vertex $w$ and the $w$-components $C_1,\dots, C_t$, respectively, are problematic if $w$ embeds into $D_T[l(u)]$ or $D_{T}[l(l(u))]$ as some $w_i$ have to be embedded into $D_T[u]$. In order to ensure that $\phi(w)$ and $\phi(w_i)$ are adjacent, we will make sure that each vertex $w_i$ is embedded in such a way that the level $L_T(\phi(w_i))$ is less than $L_T(\phi(w))+r$ in our embedding $\phi$ for some $r$ and $L_T(\phi(w))$ is small. Then the additional condition \Cref{newdef}~\ref{G4} will guarantee that $\phi(w)$ will be adjacent to $\phi(w_i)$ as desired. In addition to this, we also utilize additional types of tree partitions as supplied in \Cref{lem: feasible critical} to obtain a more sophisticated partition of $T'_0=T'[\{w\}\cup C_{[t]}]$ that is useful for obtaining a desired embedding.
In order to obtain such a partition, we crucially use the fact that the ratio between $\nu(l(u))$ and $\nu(u)$ is less than two.
These arguments are encapsulated in \Cref{strong} which ensures that we can embed a tree $T'$ into $G^2_T$ in such a way that certain vertices $x_1,x_2$ are guaranteed to be embedded so that $L_T(\phi(x_1))$ and $L_T(\phi(x_2))$ are small if some conditions are met. By exploiting this lemma, we obtain a construction of a desired universal graph with $\tfrac{14}{3}n\log_3 n + O(n) \leq 2.945 \cdot n \log_2 n$ edges as desired.

We point out an error in \cite{grahamthird} and provide a fix for it in \Cref{ch3}. Our main lemma (\Cref{strong}) for universality is proved in \Cref{sec:4}. In~\Cref{sec: construction}, we introduce a specific $(2,1)$--tree $T$ that satisfies the conditions in \Cref{strong}. Utilizing the imbalance of $\nu(l(u))$ and $\nu(u)$ in this tree and further structures of $T$, we prove in \Cref{sec: construction} that the graph $G^2_T$ has at most $\frac{14}{3}n\log_3 n + O(n)$ edges as desired. In \Cref{sec:6}, we introduce another parameter $s^{int}(n)$  as a variation of $s(n)$ and provide a lower bound for it.

\section{A remark on the proof in \cite{grahamthird}.} \label{ch3}

In this section, we discuss the mistake in \cite{grahamthird}. In order to do this, we will first briefly introduce the definitions from \cite{grahamthird}. Let {\bf $B(k)$} be a perfect binary tree of level $k$, and $V(k)$ be the vertex set of $B(k)$. Again, $B(k)$-admissible set is defined as before.

\begin{definition}\label{olddef}
For a given rooted tree $B(k)$, a graph $G(k)$ is a graph on the vertex set $V(k)$ obtained from $B(k)$ by adding all edges as follows for every $u\in V(B(k)).$

 \begin{description}
     \item[(i)] $uw$ for all descendants $w$ of $u$.
     \item[(ii)] $uw$ for all left-siblings $w$ of $u$, and $uw'$ for all descendants $w'$ of $w$.
     \item[(iii)] $uw$ for the nearest-left sibling $w$ (if it exists) of the parent of $u$, and $uw'$ for all descendants $w'$ of $w$.
 \end{description}     
\end{definition}
Note that (iii) is different from $G_{B(k)}$ as defined in \Cref{newdef}. Here, $w$ is the nearest-left sibling of the parent of $u$, instead of the nearest-left cousin. In fact, if the parent of $u$ is the left-most child of its parent, (iii) does not add any more edges in \Cref{olddef} while \Cref{newdef} does.
We say that a graph $H$ is {\bf $G(k)$-admissible} if it is an induced subgraph $G(k)[S]$ of $G(k)$ for some $B(k)$-admissible set $S$. The following is the main lemma from \cite{grahamthird}.

\begin{lemma}\label{oldlemma}
(Lemma labeled (*) in \cite{grahamthird}) 
Let $S$ be a $B(k)$-admissible set and $G=G(k)[S]$ be the subgraph of $G(k)$ induced by $S$. Let $T$ be any tree with at most $|S|$ vertices and $v\in V(T)$.
Then there exists an embedding $\lambda : T\rightarrow G$ such that the vertex set $V(G)-\lambda(V(T))$ is $B(k)$-admissible. Furthermore, the location of $\lambda(v)$ could be determined as follows:
\begin{itemize}
    \item If $|T|=|G|$, $\lambda(v)$ is root of $B(k)$.
    
    \item If $|T|<|G|$, then $\lambda(v)$ is the rightmost vertex (in terms of the DFS preorder traversal) in $\lambda(V(T))$ satisfying
    \[ |D_{B(k)}(\lambda(v))| < |T|\leq  |D_{B(k)}(\lambda(v)^*)| \]
where $\lambda(v)^*$ is the parent of $\lambda(v)$.
\end{itemize}
\end{lemma}

See \Cref{fig1} (which is same as Figure~7 in \cite{grahamthird}). It explains the part II, case (ii), subcase (b) of the proof for \Cref{oldlemma} in \cite{grahamthird}.  It was claimed in  \cite{grahamthird} that the graph $G_2$ in \Cref{fig1} is isomorphic to an $G(K)$-admissible graph, but this is not true.

Note that $l(u)$ is not a left-sibling of $u$, so $l(u)$ and $w$ are not adjacent in $G$. Thus, we cannot say that it is isomorphic to a $G(\ell)$-admissible graph. Indeed, it is easy to construct a small example that shows the claim is wrong.
For example, one can consider a $G(3)$-admissible graph on $11$ vertices with the DFS preorder traversal $x_1,\dots, x_{11}$. Then the last six vertices $\{x_6,\dots, x_{11}\}$ form a non-complete graph missing edges $x_6x_{11}, x_7x_{11}$ and $x_8x_{11}$. On the other hand the vertex sets $\{x_6, x_7, x_8\}$ and $\{x_9,x_{10},x_{11}\}$ correspond to $D_{T}[l(u)]$ and $D_{T}[u^{*}]$ as in \Cref{fig1}. However, $\{x_6,\dots, x_{11}\}$ does not induce a $G(\ell)$-admissible graph for any $\ell$ because it is routine to check that any $G(\ell)$-admissible graph having six vertices for $2\leq \ell\leq 5$ must be complete graph. 

Nonetheless, the proof can be corrected. Replace $G(k)$ with the graph $G_{B(k)}$ defined in \Cref{newdef}, then the induction process for \Cref{oldlemma} carries over without any modification. So, any $n$-vertex $G_{B(k)}$-admissible graph can be shown to be a universal graph as desired. 

In order to prove \Cref{thm1}, we only need to count the number of edges in an $n$-vertex $G_{B(k)}$-admissible graph. We use induction on $(n,k)$ to prove that any $n$-vertex $G_{B(k)}$-admissible graph contains at most $\frac{7}{2}kn + 4n$ edges. We write $(n',k') < (n,k)$ if $k'<k$ or $k'=k$ and $n'<n$.

First, as a base case, we prove a stronger bound of $\frac{7}{2}kn + n$ for the case $n=2^{k+1}-1$ with $k\in \mathbb{N}$. 
For a vertex $u$ on the level $t$ that has a left sibling, \Cref{newdef} adds at most \[\frac{n-2^{t+1}+1}{2^t} + \frac{n-2^{t}+1}{2^t} + \frac{2n-2^{t}+2 }{2^{t}} \leq \frac{4n}{2^t} \]
edges of the form $uw$ or $uw'$. On the other hand, for a vertex $u$ on the level $t$ that has no left sibling, \Cref{newdef} adds at most 
\[\frac{n-2^{t+1}+1}{2^t} + 0+ \frac{2n-2^{t}+2 }{2^{t}} \leq \frac{3n}{2^t} \]
edges of the form $uw$ or $uw'$. As there are $2^{t-1}$ vertices for each type at level $t$ for each $t\geq 1$ and the root vertex has degree $n-1$, the total number of edges in $G_{B(k)}$ is at most 
\[ (n-1) + \sum_{t\in k}(2^{t-1} \cdot \frac{4n+3n}{2^t}) \leq \frac{7}{2}kn+n, \] 
as desired.

Now we consider general $(n,k)$. Assume that the claim holds for all $(n',k')<(n,k)$, i.e., any $n'$-vertex $G_{B(k')}$-admissible graph contains at most $\frac{7}{2}k'n'+ 4n'$ edges. 
Let $G$ be an $n$-vertex $G_{B(k)}$-admissible graph, and let $v$ be the vertex with minimum level that has both children included in $V(G)$. Such a vertex $v$ exists as otherwise $G$ is a clique on at most $k+1$ vertices and it has at most $\frac{1}{2}kn\leq \frac{7}{2}kn+4n$ edges as desired. 
Let $\ell$ be the level of $v$ and $v_1,v_2$ be the left child and the right child of $v$. Letting $G_i= G[D_{B(k)}[v_i]]$, $G_1$ is isomorphic to $G_{B(k')}$ with $k'=k-\ell-1$ and $G_2$ is isomorphic to an $(n-2^{k'+1}-\ell)$-vertex $G_{B(k')}$-admissible graph. 
If $\ell>0$, then the induction hypothesis on $G-v$ yields that $G$ has at most $$(n-1)+ \frac{7}{2}(k-1)(n-1)+4(n-1) \leq \frac{7}{2}k n + 4n$$
edges as desired. Thus we may assume $\ell=0$ and $k'=k-1$. 
The edges between $V(G_1)$ and $V(G_2)$ are obtained from \Cref{newdef}~\ref{G3} and they are incident with $v_2$, children of $v_2$, children of $v_2$'s left children, children of $v_2$'s leftmost grandchildren, and so on. Thus there are at most $(2^{k}-1) + 2(2^{k}-1) + 2(2^{k-1}-1) + 2(2^{k-2}-1)+ \dots + 2(2^2-1) \leq 5(2^{k}-1)$ edges between $V(G_1)$ and $V(G_2)$.
Since $v$ has degree $n-1$, the induction hypothesis yields that $G$ contains at most 
\begin{align*}
&   (n-1) +  \left(\frac{7}{2}(k-1)(2^{k}-1)+ 2^{k}-1\right)+ \left(\frac{7}{2}(k-1)(n-2^{k}) + 4(n-2^{k}) \right) + 5 (2^{k}-1)    \\
&\leq \frac{7}{2} k n  + (4-\frac{5}{2}) n + 2(2^{k}-1) \leq \frac{7}{2} kn + 4n,
\end{align*}
as desired. Here the final inequality holds as $2^{k}\leq n$.

Finally, for any integer $n$, we take $k=\lfloor \log_2 n\rfloor$, then an $n$-vertex $G_{B(k)}$-admissible graph exists and it contains at most $\frac{7}{2}kn + 4n \leq \frac{7}{2} n \log_2 n + 4n$ edges as desired.

\section{Universality of $G^2_T$ for a $(2,1)$-balanced tree $T$.}
\label{sec:4}

In this section, we define a tree with a controlled imbalance between $\nu(u)$ and $\nu(l(u))$. Moreover, we prove that the graph $2$-generated by such an $n$-vertex tree is a universal graph  for $\mathcal{T}_n$.

\begin{definition}
    A rooted tree $T$ is called {\bf a $(K,s)$-balanced tree} (or a {\bf $(K,s)$-tree}) if it satisfies the following.
\begin{enumerate}[label=\upshape \textbf{T\arabic{enumi}}]
        \item \label{T1} 
        For each $u\in V(T)$, if $l(l(u))$ exists, then  $ \nu(l(l(u)))+\nu(l(u)) \geq \nu(u)$ holds.
        \item \label{T2}
        For each $u\in V(T)$, if $l(u)$ exists, then $ \nu(l(u)) > \frac{1}{K} \nu(u)$ holds.
        \item \label{T3} For $u,u'\in V(T)$ with $L(u)=\ell$ and $L(u')\geq \ell+s$, if $u$ is not the rightmost vertex in level $\ell$, then $\nu(u)\geq \nu(u')$.
        \item \label{T4} For any vertex $u\in V(T)$, if $u$ has a child, then $l(u)$ also has a child.
    \end{enumerate}
 \end{definition}

In other words, $T$ resembles a complete $t$-ary tree in the sense that each vertex $u$ at level $\ell$ has roughly the same  number of descendants (up to some  multiplicative constant), and a vertex on a lower level has more descendants than a vertex on a higher level. 
 Moreover, \ref{T1} ensures that if $\nu(l(u))$ is smaller than $\nu(u)$, then $\nu(l(l(u)))$ must be somewhat large, so that no two consecutively small $\nu$-values appear in level $\ell$.

Note that if $T$ is a $(K,s)$-tree, then for any $T$-admissible set $U$, $T[U]$ is also a $(K,s)$-tree. Moreover, once we list all vertices $u_1,\dots, u_t$ on the same level in $T$, we take the consecutive vertices $u_{i},\dots, u_{t}$ and let $Q=D_T[u_i]\cup \dots \cup D_T[u_{t}]$. Then the new tree $T^*$ obtained from $T[Q]$ by adding a new root $v$ is a $(K,s)$-tree.

Furthermore, if either $u_i$ and $u_{t}$ have the same parent or their parents are consecutive on the same level (i.e. $u_i^* = l(u_{t}^*)$), then the following holds for any $r\geq 0$.
\begin{equation}
\begin{minipage}{0.9\textwidth} \label{eq: tree merge}
    The induced subgraph $G^{r}_T[Q \cup \{u_{t}^*\}]$ of $G^{r}_T$ is isomorphic to $G^{r}_{T^*}$.
\end{minipage}
\end{equation}
 Indeed, it is easy to see that the \Cref{newdef}~\ref{G1}--\ref{G3} applied to $T^*$ adds exactly the same edges as in $G^r_T[Q \cup \{u_{t}^*\}]$. Furthermore, assume \Cref{newdef}~\ref{G4} adds an edge $\overrightarrow{uw}$ in $\overrightarrow{G}^r_T[Q \cup \{u_{t}^*\}]$ for the $i$-th descendant $w$ of $u'$ where $u'$ is either the $r$-th ancestor of $u$ or the nearest-left cousin of the $r$-th ancestor of $u$. Then $u'$ is still the $r$-th ancestor of $u$ or the nearest-left cousin of the $r$-th ancestor of $u$ in $T^*$ unless $L_T(u')\leq L_T(u^*_{i})$. If $L_T(u')\leq L_T(u^*_{i})$, then $u'$ is replaced by $v$. Since $v$ is the $r$-th ancestor of $u$ in $T^*$ and $w$ is $i'(\leq i)$-th descendant of $v$ in $T^*$, \ref{G4} adds $\overrightarrow{uw}$ to $\overrightarrow{G}^{r}_{T^*}$. Thus $G^r_T[Q\cup \{u^*_{t}\}]$ is a subgraph of $G^r_{T^*}$. 
 Conversely, one can similarly check that $\overrightarrow{uw} \in \overrightarrow{G}^{r}_{T^*}$ implies $\overrightarrow{uw} \in G^{r}_{T}[Q]$, hence \eqref{eq: tree merge} holds.

Moreover, again assume $u_i^*=u_t^*$ or $u_i^* = l(u_{t}^*)$. We also assume that $u_t^*$ is either the root of $T$ or the rightmost vertex on level $1$. Any $T^*$-admissible embedding $\phi:T'\rightarrow T^*$ whose image $\phi(T')$ does not contain the root $v$  satisfies that $V(T)-\phi(T')$ is either a $T$-admissible set or a $T$-admissible set plus a vertex $u^*_t$, where the latter holds only if the $T^*$-admissible set has at least $\nu_T(u^*_t)$ vertices.
Thus, when we want to find a $T$-admissible embedding of a tree $T'$ into $G_T$ with $|T'|\geq\nu_T(u_t^*)$, we can instead find a $T^*$-admissible embedding $\phi$ of $T'-x$ into $G^r_{T^*}$ for some vertex $x$ and assign $\phi(x)=u_t^*$. Because $|T'|\geq \nu_T(u_t^*)$ and $u_t^*$ is adjacent to all vertices in $G^{r}_{T}[Q]$, this yields a $T$-admissible embedding of $T'$ into $G^{r}_T$. In the proof of the following lemma, we will frequently use this trick to find a $T$-admissible embedding.

\begin{lemma}\label{strong}
    Let $T$ be an $n$-vertex $(2,1)$-balanced tree rooted at $v$ with the level $L(T)=k$ and $T'$ is a tree on $n'\leq n$ vertices.  Let $v_1,\dots, v_{t}$ be the children of $v$ from left to right.      For two (not necessarily distinct) vertices $x_1,x_2 \in V(T')$, there exists a $T$-admissible embedding $\phi : T'\rightarrow G^2_T$ satisfying the following.
    \begin{enumerate}[label=\upshape \textbf{$\Phi$\arabic{enumi}}]
        \item $\phi(x_1)$ has the minimum level among the vertices in $\phi(T')$.
        \label{strong1}
         \item 
         If $t=2$ and $n-2\geq n'\geq \nu(v_t)\geq 2$, then $L_T(\phi(x_2))\leq 2$. \label{strong2}
    \end{enumerate}      
\end{lemma}
\begin{proof}
We write ${\bf A_{k}}$ to denote the statement that the lemma holds for all trees $T$ with level at most $k$, and we write ${\bf A_{k,t}}$ to denote the statement that the lemma holds for all trees with level at most $k$ whose root has at most $t$ children.
We write $(k',t')< (k,t)$ if $k'<k$ or $k'=k$ and $t'<t$.

We use induction on $(k,t)$. When $k \leq  2$, then $G^2_T$ is a complete graph, so it is easy to check that ${\bf A_k}$ holds. Assume that $k\geq 3$ and ${\bf A_{k',t'}}$ holds for all $(k',t')< (k,t)$. \newline 

First, if $\nu(v_t)=1$, then we use ${\bf A_{k,t-1}}$ to find a $T$-admissible embedding $\phi$ of $T-x_1$ into $G^2_T -v_t$.
If $n'<n$, we extend it by embedding $x_1$ to $v_t$, then $\phi$ is a desired embedding. If $n'=n$, take $u'=\phi^{-1}(v)$ and we modify $\phi$ by letting $\phi(x_1)=v$ and $\phi(u')=v_t$. As both $v$ and $v_t$ are adjacent to every vertex in $G^2_T$, this yields a desired $T$-admissible embedding satisfying \ref{strong1}. In any of the cases, \ref{strong2} is vacuous. Thus $\phi$ is a desired embedding. Therefore, from now on, we assume $\nu(v_t)\geq 2$. This with \ref{T4} further implies that $\nu(v_i)\geq 2$ for all $i\in [t]$. 

If $n' < \nu(v_t)$, then we use ${\bf A_{k-1}}$ to find a $T$-admissible embedding of $T'$ into $D_T[v_t]$ satisfying \ref{strong1}. As \ref{strong2} is vacuous in this case, this yields a desired embedding. Thus, we assume $n'\geq \nu(v_t)$ in each of the following cases.  \newline

For the cases of $t\in \{1,2\}$, we delete either the root $v$ or two vertices $v$ and $v_1$ to obtain a tree with one less level. As the remaining tree has a smaller level, we can invoke the induction hypothesis to embed $T'$, with the possible exception of a few vertices. We can manually embed these a few unembedded vertices in such a way that \ref{strong1} and \ref{strong2} are satisfied. \newline 


\noindent {\bf Case 1. ${\bf t=1}$.} We consider the tree $D_{T}[v_1]$ of level $k-1$. If  $n'\leq n-1$, we use ${\bf A_{k-1}}$ to find a $T$-admissible embedding $\phi$ of $T'$ into $D_{T}[v_1]$ so that \ref{strong1} holds. As $t=1$, \ref{strong2} is vacuous, thus this yields a desired embedding. 

If $n'=n$, then  use ${\bf A_{k-1}}$ to find  a $T$-admissible embedding $\phi$ of $T'-x_1$ into $D_{T}[v_1]$.  We let $\phi(x_1)=v$, then we obtain a desired $T$-admissible embedding $\phi$ satisfying \ref{strong1} and \ref{strong2} is vacuous. \newline

\noindent {\bf Case 2. ${\bf t=2}$.} Let $u_1,\dots, u_a$ be the children of $v_1$ from left to right and let $u_{a+1},\dots, u_{a+b}$ be the children of $v_2$ from left to right. We know that $a,b\geq 1$ as $\nu(v_1), \nu(v_2)\geq 2$.  We further divide the cases according to the value of $n'$. Recall that we have assumed $n'\geq \nu(v_2)$.
\newline

\noindent {\bf Case 2-1. ${\bf t=2}$ and ${\bf \nu(v_2)\leq  n'\leq  n-2}$.} In this case, we consider a new tree $T^*$ obtained from $T$ by removing $v,v_1,v_2$ and adding a new root $v'$ that has $u_1,\dots, u_{a+b}$ as children from left to right. By \eqref{eq: tree merge}, $G^2_{T^*}$ is isomorphic to the subgraph $G^2_T-\{v,v_1\}$ of $G^2_T$.
As $T^*$ has the level $k-1$, we use ${\bf A_{k-1}}$ to obtain a $T^*$-admissible embedding $\phi$ of $T'-x_1$ into $G^2_T-\{v,v_1\}$ such that the level of $\phi(x_2)$ is the minimum among $\phi(T'-x_1)$. As $n'\geq \nu(v_2)$ and  $\phi$ is $T^*$-admissible, $\phi(x_2)$ is on level one in $T^*$, thus $\phi(x_2)\in \{u_1,\dots, u_{a+b}\}$. 
 Moreover, as $|T'|-1= n'-1 \leq n-3$, the vertex $v_2$  (which corresponds to the root of $T^*$) is not in $\phi(T'-x_1)$. By letting $\phi(x_1)=v_2$, the map $\phi$ extends to an embedding of $T'$ into $G^2_T$, yielding  a desired $T$-admissible embedding satisfying \ref{strong1} and \ref{strong2}. Note that this is a $T$-admissible embedding because $n'\geq \nu(v_2)$.
 \newline

\noindent {\bf Case 2-2. ${\bf t=2}$ and ${\bf n'\in \{n-1,n\}}$}.
Again, we consider the new tree $T^*$ obtained from $T$ by removing $v,v_1,v_2$ and adding a new root $v'$ that has $u_1,\dots, u_{a+b}$ as children from left to right. 
Take a leaf $w$ (or an isolated vertex if there's no leaf) of $T-x_1$ and let $w'$ be the unique neighbor of $w$ in $T-x_1$. If $w$ is an isolated vertex in $T-x_1$, then let $w'$ be an arbitrary vertex. 
By \eqref{eq: tree merge}, $G^2_{T^*}$ is isomorphic to the subgraph $G^2_T-\{v,v_1\}$ of $G^2_T$.
Apply the induction hypothesis ${\bf A_{k-1}}$ to find a $T^*$-admissible embedding $\phi$ of $T-\{x_1,w\}$ into $G^2_{T}-\{v,v_1\}$ satisfying \ref{strong1}, where $w'$ plays the roles of $x_1$, respectively. Then we know that $L_T(\phi(w'))=L_{T^*}(\phi(w'))+1 \leq 2$.
We modify this embedding by letting 
\[\phi(w) = v_1 \enspace \text{ and } \enspace \phi(x_1)=\left\{\begin{array}{ll}
    v_2 & \text{ if } n'=n-1. \\
    v & \text{ if }n'=n.
\end{array}\right.  \]
Then $\phi(w)$ and $\phi(w')$ are adjacent by \Cref{newdef} because $\phi(w)=v_1$ is either the parent of $\phi(w')$, the left sibling of $\phi(w')$, or the nearest-left cousin of the parent of $\phi(w')$. Also, \ref{strong2} is vacuous in this case. Thus, this yields a desired $T$-admissible embedding. \newline

\noindent {\bf Case 3. ${\bf t\geq 3.}$} In this case, we do not care whether our embedding satisfies \ref{strong2} as it is vacuous.
Let 
\[ x= \nu(v_t), \enspace y=\nu(v_{t-1}) \enspace \text{and} \enspace z=\nu(v_{t-2}).\]
Recall that $x,y,z\geq 2$. Based on the values of these numbers, we take a vertex $w$ of $T$ and collection of $w$-components satisfying certain conditions as the following two cases.
\begin{description}[leftmargin=3em, style=nextline]
  \item[{\bf Case 3-1.}] If there exists a $(x,y)$-feasible collection $C_1,\dots, C_p$ of $w$-components for some vertex $w$ with $x_1\notin C_{[p]}$, then we choose it.
   \item[{\bf Case 3-2.}] If such a collection does not exist, then we choose a $(x,y)$-critical collection $C_1,\dots, C_p$ of $w$-components for some $w$ with $x_1\notin C_{[p]}$.
\end{description}
By \Cref{lem: feasible critical}, one of the above cases must occur. Let $w_i$ be the unique neighbor of $w$ in $C_i$ for each $i\in [p]$. 
Let 
\[T'_0 = T'[ \{w\} \cup C_{[p]}] \enspace \text{ and }\enspace T'_1= T'-T'_0.\]
Note that all the edges between $T'_0$ and $T'_1$ are incident to $w$. \newline

We proceed differently in the two cases. 
Essentially, our plan is to find a vertex $w$ and a collection $C_1,\dots, C_p$ of $w$-components of an appropriate size so that $T'_0=T'[C_{[p]}\cup \{w\}]$ fits inside $D_{T}[v_{t-1}]\cup D_{T}[v_t]$. Once we obtain such a subtree $T'_0$, then we invoke induction hypothesis to embed $T'_0$ into $D_{T}[v_{t-1}]\cup D_{T}[v_t]$ in such a way that $w$ embeds to $v_{t}$. This is Case 3-1. 
If it is not possible to find such $(w,C_{[p]})$, then we are in Case 3-2 and we instead find such a tuple with $|C_{[p]}|$ slightly bigger than what we want. As $|C_{[p]}|$ is bigger, we can further find a vertex $w'\in C_1$ and $w'$-components $C'_1,\dots, C'_{p'}$ of $T'_0[C_1\cup \{w\}]$ with appropriate sizes. With this choice, $C_{[p]}-C'_{[p']}$ will fit into $D_{T}[v_{t}]\cup D_{T}[v_{t-1}]$ and the remaining $C'_{[p']}\cup \{w\}$ will fit into $D_{T}[v_{t-1}]\cup D_{T}[v_{t-2}]$. See \Cref{fig:Case3-2}. We will carefully embed $w$ and $w'$ so that the edges between those subgraphs of $T'_0$ will also map to the edges in $G^2_T$. 
In particular, $C_{1}-C'_{[p']}$ contains two special vertices $w'$ and $w_1$ sending edges outside, so we will use \ref{strong2} to control on where to embed them. For this, we will keep track of the number of vertices remaining after embedding each piece into $D_{T}[v_t]$ to make sure the conditions to invoke \ref{strong2} are satisfied.
\newline

\noindent {\bf Case 3-1. ${\bf (w,C_1,\dots, C_p)}$ is ${\bf (x,y)}$-feasible.} 
In this case, we know \(x\leq |T'_0|\leq x+y-2.\)
Take $D_T[v_{t-1}]\cup D_T[v_t]$ and add a new root $v'$ as the parent of $v_{t-1}$ and $v_{t}$ to obtain a new tree $T^*$. 
For $Q=D_T[v_{t-1}]\cup D_T[v_t]\cup \{v\}$, 
\eqref{eq: tree merge} ensures that $G^2_T[Q]$ is isomorphic to $G^2_{T^*}$. Hence we can use the induction hypothesis ${\bf A_{k,2}}$ with $w$  playing the role of $x_1$. This yields a $T^*$-admissible embedding $\phi_0$ of $T'_0$. 
As $|T'_0|\leq x+y-2 < |Q|$, we know $v'\notin \phi_0(T'_0)$. Thus, $\phi_0$ is actually a $T$-admissible embedding of $T'_0$ into $G^2_T$ as well.  Furthermore, $\phi(w)$ is on level $1$, hence $\phi(w)=  v_t$. 

As $|T'_0|\leq x+y-2$, the tree $T_1=T\setminus\phi(T'_0)$ is a subtree of $T$ with the root having $t-1$ children.
Hence, we can apply the induction hypothesis ${\bf A_{k,t-1}}$. 
This yields  a $T_1$-admissible embedding $\phi_1$ of $T'_1$ into $G^2_T-\phi_0(T'_0)$, where $\phi_1(x_1)$ has the minimum level. 
Because all the edges between $T'_0$ and $T'_1$ are incident with $w$ and $\phi_0(w) = v_t$ is adjacent to all vertices in $\phi_1(T'_1)$, combining $\phi_0$ and $\phi_1$ yields a $T$-admissible embedding $\phi$. 

This is our desired embedding except one case, when $\phi(x_1)$ has the level one and it is not the minimum. It happens only when $w=x_1$ and $n'=n$ hold. 
If this happens, take a vertex $u'=\phi^{-1}(v)$ which maps to the root $v$. As we know $\phi(w)=v_t$, we modify $\phi$ so that $\phi(x_1)=v$ and $\phi(u')= v_{t}$. As both $v$ and $v_{t}$ are adjacent to all vertices in $G^2_T$, this modification also yields an $T$-admissible embedding.

The final $T$-embedding $\phi$ satisfies \ref{strong1}. Indeed, $x_1$ is either $w$ or in $T'_1$. In either case, it has level one if $n'<n$ and zero if $n'=n$.   \newline

\noindent{\bf Case 3-2. ${\bf (w,C_1,\dots, C_p)}$ is ${\bf (x,y)}$-critical.} 
In this case, we further partition $T'_0$ into several forests.
Note that $y<x$ holds in this case as otherwise application of \Cref{lem2} with $x-1$ playing the role of $x$ yields a ${(x,y)}$-feasible collection, meaning that we would be in Case~3-1.  

We define set $C_0$ as follows.
\[ 
C_0 = \left\{ \begin{array}{ll}
\{w\} & \text{ if}\enspace |T'|\geq x+y+z,\\
  T'- C_{[p]}  &  \text{ if}\enspace |T'|\leq  x+y+z-1,\\
\end{array}\right.
\]
In either case, $w$ is in $C_0$. 

\begin{figure}[h]
\centering

\includegraphics[width=0.6\textwidth]{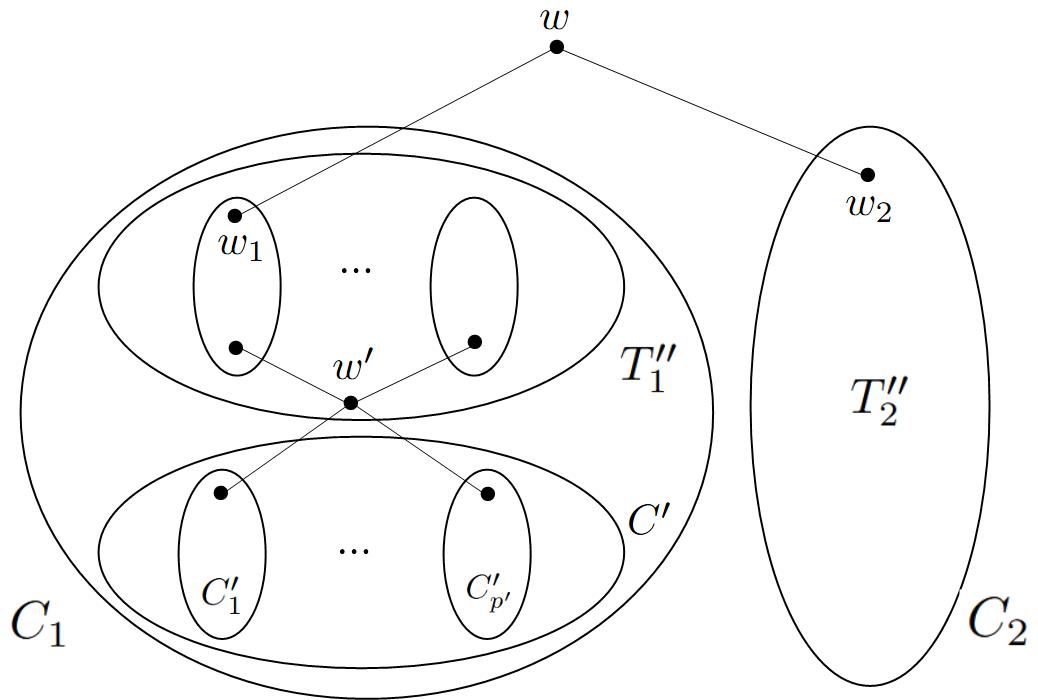}
\caption{Structures of $T'_{0}$}\label{fig:Case3-2}
\end{figure}

Recall that $(x,y)$-criticality implies $|C_i|\geq y$, thus  we have $|C_{[2]}|\geq 2y > x$ while $|C_{I}|\leq x-2$ for all $I\subsetneq [p]$.
This yields that we have $p=2$. Thus, $2y\leq |C_1|+|C_2| \leq |C_{[2]}|\leq 2x-3$ implying $y\leq x-2$. 
Moreover, $(x,y)$-criticality implies the following.
\[x+y-1 \leq |T'_0| \leq 2x-2 .\]
Note that $y\leq x-2$ implies $x+y-1< 2x-2$. Let $a= \min\{ 1, |T'_0|-x-y+1\}$. In particular,  we have $a=0$ if $|T'_0|=x+y-1$ and $a=1$ if $|T'_0|=2x-2>x+y-1$, hence $|T'_0|\leq 2x-3+a$.
From this, we have the following. 
\[1\leq |T'_0|-(x+y)+2-a \leq  x-y < y.\] 
Because $|C_1\cup \{w\}|\geq y+1$ we can 
apply \Cref{lem2} to $T'_0[C_1\cup \{w\}]$  
to obtain a vertex $w'$ and $w'$-components $C'_1,\dots, C'_{p'}$ of $T'_0[C_1\cup \{w\}]$ with $w\notin C'_{[p']}$
satisfying the following where $C'=C'_{[p']}.$
\begin{align}\label{eq: C'}
    |T'_0|-(x+y)+2-a\leq \left|C'\right| \leq  2|T'_0|-2(x+y)-2a+3.
\end{align}
Because $|T'_0|\leq 2x-3+a$ and $x\leq 2y-1$ hold, \eqref{eq: C'} implies 
\begin{align}\label{eq: T'0-C' sizes}
x  \leq  2y -1 \leq   2(x+y)-|T'_0|+2a-4 \leq  |T'_0 - \left(C'\cup \{w\}\right)| \leq x+y-2.
    \end{align}
In particular, we know $w'\neq w$. Indeed, if $w'= w$, then $C'_{[p']}$ is $C_1$ , thus $|T'_0 - \left(C'\cup \{w\}\right)|=|C_2| \leq x-1$, which contradicts the 
$(x,y)$-criticality.

Let 
\[T''_0= T'[C'\cup C_0], \enspace T''_1= T'_0[C_1]-C' \enspace \text{ and }\enspace T''_2 = T'[C_2].\]
In particular $w'$ is in $T''_1$. As $w_1$ and $w$ are adjacent and $C'$ does not contain $w$, we know that either $w'=w_1$ or $w_1$ is in the same $w'$-component with $w$ in $C_1\cup \{w\}$. In any case, we have $w_1\in V(T''_1)$.
Furthermore, as $T''_1\cup T''_2 = T'_0 - \left(C'\cup \{w\}\right)$, \eqref{eq: T'0-C' sizes} implies 
\begin{align}\label{eq: T'' sizes}
 x\leq |T''_1\cup T''_2| \leq x+y-2
\end{align}

\begin{figure}[h]
\centering

\includegraphics[width=0.6\textwidth]{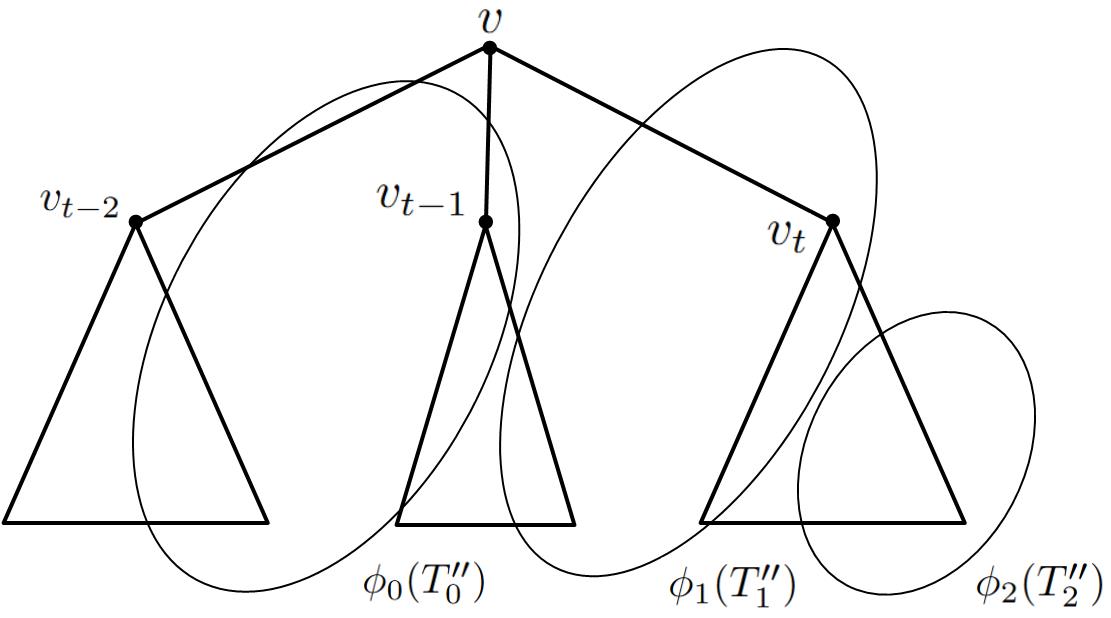}
\caption{Each tree $T''_i$ embeds as above. In addition, $v_t= \phi_1(w')$ and $v_{t-1}\in \{\phi_0(x_1), \phi_0(w)\}$.}
\end{figure}

 We plan to embed $T''_2, T''_1, T''_0$ and $T'_1-C_0$ in order. First, we embed $T''_2$.  Let $u'_1,\dots, u'_{a}$ be the children of $v_t$ and consider the tree $D_T[v_t]$.
As $T$ is $(2,1)$-tree and $\nu(v_{t-1})=y$, \ref{T3} ensures $\nu(u'_i)\leq y$ for all $i\in [a]$ and the definition of $(x,y)$-critical collection implies $|T''_2|\geq y$. Hence, we can apply the induction hypothesis ${\bf A_{k-1}}$ with $w_2$  playing the role of $x_1$ to embed $T''_2$ into $D_T[v_t]$ in such a way that $w_2$ is embedded into a child of $v_t$. We denote such embedding $\phi_2$. Indeed, \ref{strong1} ensures that the image of $w_2$ has the minimum level, thus $T$-admissibility and the fact $|T''_2|\geq \nu(u'_a)$ ensures this.  \newline

Now we embed $T''_1$. 
Let the tree $T^*$ be obtained from $D_T[v_{t-1}]\cup D_T[v_{t}] -\phi_2(T''_{2})$ by adding the root $v'$.
Then the definition of $(x,y)$-critical collection implies $\left|T''_2\right| \leq  x-2 = \nu_T(x_t)-2$, hence $\nu_{T^*}(v_t)\geq 2$.
Let $Q_1= \{v\}\cup D_T[v_{t-1}]\cup D_T[v_{t}] -\phi_2(T''_{2})$, then  \eqref{eq: tree merge} implies that $G^{2}_T[Q_1]$ is isomorphic to $G^{2}_{T^*}$.  Furthermore, 
 \eqref{eq: T'' sizes} implies  
 \[\nu_{T^*}(v_{t}) = x - |\phi_2(T''_2)|\leq  |T''_1\cup T''_2|-|T''_2|= |T''_1| \leq x+y-2 - |\phi_2(T''_2)| = \nu_{T^*}(v_t)+ \nu_{T^*}(v_{t-1})-2.\] Hence, we have $|T^*|-2\geq |T''_1|\geq \nu_{T^*}(v_t)\geq 2$, so  we can apply the induction hypothesis ${\bf A_{k,2}}$ with $w', w_1$ playing the roles of $x_1,x_2$, respectively, to find a $T$-admissible embedding $\phi_1$ into $G^2_T-\phi_2(T''_2)$
 such that $\phi_1(w')=v_{t}$ and $L_T(\phi_1(w_1))\leq 2$.  Note that we use \ref{strong2} here to ensure both vertices $w'$ and $w_1$ are mapped to vertices of small level. \newline 
 
Now we embed $T''_0$. Let $T^{**}$ be the tree obtained from $D_T[v_{t-2}]\cup D_{T}[v_{t-1}]\setminus (\phi_1(T''_1)\cup \phi_2(T''_2))$ by adding the root $v'$. Since $|T''_1| \leq \nu_{T^*}(v_t)+\nu_{T^*}(v_{t-1})-2$, we have $\nu_{T^{**}}(v_{t-1}))\geq 2$.
Let $Q_0= \{v\}\cup D_T[v_{t-2}]\cup D_T[v_{t-1}]\setminus (\phi_1(T''_1)\cup \phi_2(T''_2))$, then  \eqref{eq: tree merge} implies that $G^{2}_T[Q_0]$ is isomorphic to $G^{2}_{T^{**}}$.  
Furthermore, by our choice of $C_0$, we know $|T''_0|+|T''_1|+|T''_2|\leq |T'_0|+|C_0|-1 \leq \max\{ 2x-2, x+y+z-1\} = x+y+z-1$ because $x\leq y+z$.
Thus $|T^{**}|-2\geq |T''_0| \geq \nu_{T^{**}}(v_{t-1})\geq 2$.
Therefore, we can apply the induction hypothesis ${\bf A_{k,2}}$ according to the following two subcases.

\begin{enumerate}
    \item{}
If $x_1\in C_0-\{w\}$, then we apply the induction hypothesis ${\bf A_{k,2}}$ with $x_1, w$ playing the roles of $x_1, x_2$ to find a $T$-admissible embedding $\phi_0$ of $T''_0$ into $G^2_T[Q_0]$ such that $\phi_0(x_1)=v_{t-1}$ and $L_T(\phi_0(w))\leq 2$. 
    \item{}
If $x_1\notin C_0-\{w\}$, then either $T'_1-C_0$ contains $x_1$ or $w=x_1$. In either case, we apply the induction hypothesis ${\bf A_{k,2}}$ with $w$ playing the role of $x_1$ to find a $T$-admissible embedding $\phi_0$ of $T''_0$ into $G^2_T[Q_0]$ such that $\phi_0(w)=v_{t-1}$. 
\end{enumerate}

Finally, we embed $T'_1-C_0$. If $T'_1-C_0$ is empty, we let $\phi'$ be the empty embedding. 
Otherwise, let $T^{\#}$ be the tree $T - \bigcup_{i=0}^{2} \phi_{i}(T''_i)$. 
We apply the induction hypothesis ${\bf A_{k,t-2}}$ to obtain a desired $T$-admissible embedding $\phi'$ of $T'_1-C_0$ into $G_{T^{\#}}^2 = G^2_T - \bigcup_{0\leq i\leq 2}\phi_i(T''_i)$, such that $\phi'(x_1)$  has the minimum level if $x_1\in T'_1-C_0$.  \newline

Let the final embedding $\phi$ be the embedding obtained by concatenating $\phi_2, \phi_1, \phi_0$ and $\phi'$. Note that all the edges between $T''_0, T''_1, T''_2$ and $T'_1-C_0$ are incident to one of two vertices $w, w'$. As $\phi(w')=v_{t}$ is adjacent to all vertices in $G^2_T$ so those edges incident with $w'$ map to an edge in $G^2_T$.

In any case, $\phi(w)$, $\phi(w_1)$ and $\phi(w_2)$ also have level at most two. Thus they are all children of a vertex in $\{v_{t-2}, v_{t-1}, v_t\}$, they share the common grandparent and they are all adjacent. This shows that $\phi$ is indeed a $T$-admissible embedding. 

Finally, we check \ref{strong1} for our embedding. We know that either $x_1\in T'_1-C_0$ or $x_1\in C_0-\{w\}$ or $x_1=w$.
\begin{enumerate}
    \item If $x_1\in T'_1-C_0$, we obtain \ref{strong1} from the choice of $\phi'$.
    \item  If $x_1\in C_0-\{w\}$, it means $C_0-\{w\}$ is nonempty, so $|T'|\leq x+y+z-1$. Since $v\notin \phi(T')$ and that $x_1$ is embedded at level $1$, we also obtain \ref{strong1}. 
    \item If $x_1=w$ and $n'<n$,  the level of $\phi(w)=v_{t-1}$ is one, so we obtain \ref{strong1}.
    \item If $x_1=w$ and $n'=n$, consider the vertex $u'=\phi^{-1}(v)$ that maps to the root $v$. The vertex $u'$ belongs to $T'_1$, thus it has neighbors only in $T'_1\cup \{w\}$ and $\phi(T'_1\cup \{w\})\subseteq N_{G^2_T}(v_{t-1})\cup \{v_{t-1}\}$ by our choice of $\phi'$. Thus we swap the values of $\phi$ so that $\phi(x_1)=v$ and $\phi(u')=v_{t-1}$, then it is still a $T$-admissible embedding of $T'$ into $G^2_T$, and it additionally satisfies \ref{strong1}.
\end{enumerate}
 This proves ${\bf A_{k,t}}$ and the lemma.
\end{proof}

\section{A construction of universal graphs}\label{sec: construction}

We now construct a $(2,1)$-tree $T$ with level $L(T)=k$ where $G^2_T$ contains at most $\frac{14}{3\log 3}n\log n+O(n)$ edges. This construction with \Cref{strong} proves \Cref{thm2}. 
We construct this tree  level by level while assigning `type' to each vertex.

\begin{figure}[h]
\centering

\includegraphics[width=1.0\textwidth]{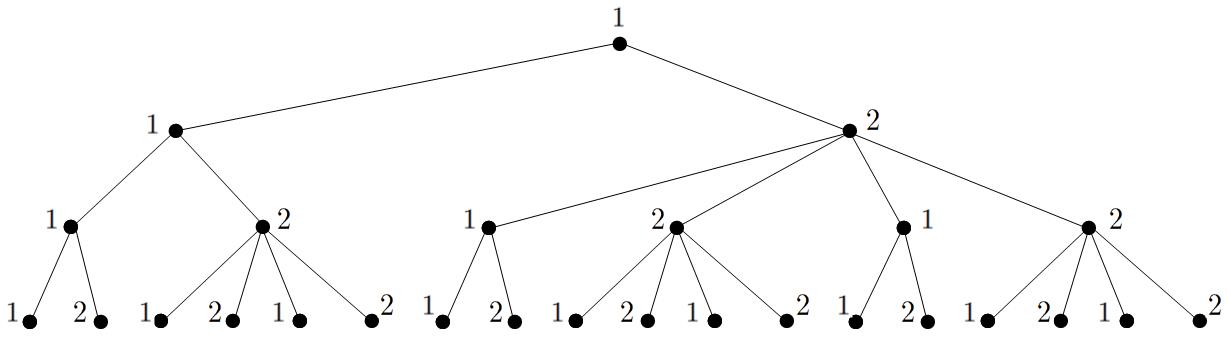}
\caption{A few levels of $T_{k}$}
\end{figure}

We define a tree $T_{k}$ of level $k$ as follows.  We start with a root vertex $v$ of type $1$. For every vertex $x$ on level $\ell<k$ of type $p\in\{1, 2\}$, we add vertices as its children as follows.
\begin{enumerate}
    \item[$\bullet$] If $p=1$, add exactly two children of type $1, 2$ from left to right, respectively.
    \item[$\bullet$] If $p=2$, add exactly four children of type $1, 2, 1, 2$ from left to right, respectively.
\end{enumerate}
Repeating this until level $k$, then we obtain our tree $T_{k}.$
Let $x^{\ell}_1,\dots, x^{\ell}_{t_{\ell}}$ be the vertices on level $\ell$ ordered from left to right. 

Note that this definition ensures that all vertices of the same type on the same level have the same number of descendants. 
We define the sequence $a_{\ell,p}=a^{k}_{\ell,p}$ be the number $\nu_{T_{k}}(x)-1$ of descendants of a vertex $x$ at level $\ell$ of type $p$. We will omit the superscript $k$ as it is clear from the context. Then the above definition yields the following recursive relation on the sequence $a_{\ell,p}$.
\begin{align*}
a_{\ell, p}=
\begin{cases}
0 & \text{ if } \ell=k \\
a_{\ell+1, 1}+a_{\ell+1, 2}+2 & \text{ if } \ell <k, p=1\\
2a_{\ell+1, 1}+2a_{\ell+1, 2}+4 & \text{ if } \ell <k, p=2\end{cases}
\end{align*}
From this relation, it is routine to check the following and it is obvious that $T_k$ is indeed a $(2,1)$-tree from these.
\begin{enumerate}[label=\upshape \textbf{R\arabic{enumi}}]
    \item $a_{\ell, 2}=2a_{\ell, 1} = 2 (3^{k-\ell}-1)$.
    \label{R1}
    \item For $\ell\in [0,k]$, the types for the vertices $x^{\ell}_1,\dots, x^{\ell}_{t_{\ell}}$ are $1,2,\dots, 1,2$, respectively. In particular, if $\ell>0$, then  the number of type-$1$ vertices and type-$2$ vertices are equal.
    \label{R2}
\end{enumerate}

As we have shown that $T_k$ is a $(2,1)$-tree, it only remains to count the number of edges in $G_{T_k}^{2}$. We instead count directed edges in the digraph $\overrightarrow{G} =\overrightarrow{G}_{T_{k}}$ which will provide an upper bound on $e(G_{T_k})$. We say a directed edge is an edge of type $i$ if it is added from \Cref{newdef}~${\bf Gi}$. As we take $\overrightarrow{G}_{T_{k}}$, there are no edges of type $4$, so all edges are type $1,2$ or $3$. We will later add the number of edges of type $4$ to $e(\overrightarrow{G}_{T_k})$, to count $e(\overrightarrow{G}^{2}_{T_k})$.

For two consecutive vertices $x^{\ell}_{2t-1}, x^{\ell}_{2t}$ in level $\ell\geq 1$ with $t\in [t_{\ell}]$, we count the number of directed edges  from their children. First, we count the edges of type $1$ and $2$. The six children of $x^{\ell}_{2t-1}, x^{\ell}_{2t}$ send directed edges towards their own descendants and descendants of their left siblings. Thus the vertices from left to right send
$$a_{\ell+1,1},\enspace  a_{\ell+1,1}+a_{\ell+1,2} +1 ,\enspace a_{\ell+1,1},\enspace  a_{\ell+1,1}+a_{\ell+1,2}+1,\enspace 2 a_{\ell+1,1} + a_{\ell+1,2}+2,\enspace 2a_{\ell+1,1}+ 2 a_{\ell+1,2}+3$$ edges of type $1$ or $2$, respectively.
This sums up to  

   \[ 8 a_{\ell+1,1}+ 5 a_{\ell+1,2}+ 7 =18 a_{\ell+1,1} + 7\leq 2(\nu(x^{\ell}_{2t-1})+\nu(x^{\ell}_{2t}))+ 7.\]
Hence, the total number of edges of type $1$ or $2$, coming from vertices on level $\ell+1$ is at most
\begin{align}\label{eq: type12}
2\sum_{i\in [t_{\ell}]} \nu(x^{\ell}_i) + \frac{7}{2} t_{\ell} \leq 2n + 4 t_{\ell}.
\end{align}
Here, the final inequality holds because $\sum_{i\in [t_{\ell}]}\nu(x^{\ell}_i)\leq n$. 

Second, we count the edges of type $3$. Each child of $x^{\ell}_{2t-1}$ sends $\nu(x^{\ell}_{2t-2})$ edges of type $3$, and each child of $x^{\ell}_{2t}$ sends $\nu(x^{\ell}_{2t-1})$ edges of type $3$. Hence, the number of edges outgoing from the vertices on level $\ell+1$ is at most 
$$\sum_{i\in [t_{\ell}/2]} \left(2\nu(x^{\ell}_{2i-2})+ 4\nu(x^{\ell}_{2i-1})\right) \leq
\frac{8}{3}\sum_{i\in [t_{\ell}/2]} \left(\nu(x^{\ell}_{2i-1})+ \nu(x^{\ell}_{2i})+\frac{3}{8}\right) \leq \frac{8}{3} n+t_{\ell}.$$
Here, the first inequality holds because $\nu(x^{\ell}_{2i-2})=\nu(x^{\ell}_{2i})$ and $\tfrac{2}{3}\nu(x^{\ell}_{2i})+1\geq\tfrac{4}{3}\nu(x^{\ell}_{2i-1})$ from \ref{R1} and the final inequality holds because  $\sum_{i\in [t_{\ell}]}\nu(x^{\ell}_i)\leq n$. 
This together with \eqref{eq: type12} yields that at most $\frac{14}{3} n + 5t_{\ell}$ edges of type $1,2,3$ come from the vertices in level $\ell+1$. 

As $n=|T_k| = a_{0,1}+1 = 3^k$,  we have $k=\log_{3}{n}$. The three vertices on levels of at most one send at most $3n$ edges, we have the following in total
\begin{align}
    e(\overrightarrow{G}_{T_{k}}) \leq 3 n + \frac{14}{3}kn + 5 \sum_{\ell=1}^{k-1} t_{\ell}  \leq \frac{14}{3} n \log_3 n + 8 n.
    \label{eq: edges in Trk}
\end{align}
Here, we obtain the final inequality since $\sum_{\ell=0}^{k-1} t_{\ell}\leq n$.  \newline

We now consider general $n$-vertex $T_k$-admissible subgraphs of $G_{T_k}$.
For a given sufficiently large $n$, consider an integer $k$ satisfying
$$3^{k-1}<n\leq 3^k.$$
Let $U$ be the $n$-vertex $T_{k}$-admissible set in $V(T_{k})$ and let $G$ be the subgraph of $G_{T_{k}}$ induced by $U$.
Let $X\subseteq U$ be the set of vertices $x$ such that $D_T[x] \subseteq U$ while $D_T[x^*]\not\subseteq U$ for its parent $x^*$.

The following claim will be useful to bound the size of $U\setminus D[X]$ where we write $D[X]=\bigcup_{x\in X} D_T[x]$.
\begin{claim}
The set $U\setminus D[X]$ contains at most one element from each level. 
\end{claim}
\begin{proof}
Consider the DFS preorder traversal $x_1,\dots, x_{3^k}$ of $T_{k}$. As $U$ is a $T_{k}$-admissible set, $U= \{x_1,\dots, x_{n}\}$. Assume that $x_{i}, x_{j} \notin D[X]$ with $i<j\leq n$ and they have the same level, then all vertices $x_{i+1},\dots, x_{j-1}$ must also be in $U$. However, by the definition of the DFS preorder traversal, DFS visits $x_j$ only after it visits all descendants of $x_i$. This implies $D_{T}[x_i]\subseteq U$, hence $x_i \in D[X]$, a contradiction.  This proves the claim.
\end{proof}

We again count the number of directed edges in $\overrightarrow{G}=\overrightarrow{G}_{T_{k}}[U]$. Note that this does not include the edges of type $4$, so we will add such edges at the end.
For each vertex $u\in U\setminus D[X]$ of level $\ell$, it sends at most $\nu(u^*)+ \nu(l(u^*))$ outgoing edges towards the vertices in $U$ where $u^*$ is the parent of $u$. 
This is at most $2a_{\ell-1,2}+2 \leq 3a_{\ell-1,2}$. As there is at most one element of level $\ell$ in $U\setminus D[X]$, the number of directed edges coming from $U\setminus D[X]$ is at most 
\begin{align}\label{eq: edges in U1}
    (n-1)+\sum_{\ell\in [k]} 3a_{\ell-1,2} \leq n + \sum_{\ell\in [k]} 2\cdot 3^{k-\ell+2} \leq 28n.
\end{align}
Here, the first term $n-1$ counts the number of edges from level $0$ and the final inequality holds as $n\geq 3^{k-1}$.

Now we count the number of edges within $D[X]$. Consider a vertex $x\in X$ of level $\ell$. We count the number of directed edges within $D_x=D_T[x]$. 


\begin{description}
    \item[Case 1: $x$ is of type $1$:]
 For $D_x=D_T[x]$, the subgraph induced by $D_x$ is isomorphic to $T_{k-\ell}$. Thus \eqref{eq: edges in Trk} yields that  the number of edges within $D_x$ is at most 
\[e(\overrightarrow{G}[D_x]) \leq \frac{14}{3}|D_x|\log_3 |D_x| + 8|D_x|.\]
    \item[Case 2: $x$ is of type $2$:]

Denote the children of $x$ as $x'_1, x'_2, x'_3, x'_4$ from left to right. Let $G_1$ be a subgraph induced by $\{x\}\cup D_{T}[x'_1]\cup D_{T}[x'_2]$, and $G_2$ be a subgraph induced by $\{x\}\cup D_{T}[x'_3]\cup D_{T}[x'_4]$. The graphs $G_1, G_2$ are both isomorphic to $T_{k-\ell}$, so the number of edges within $G_i$ is at most $\frac{14}{3}|G_i|\log_3 |G_i| + 8|G_i|.$

Now we count the edges between $G_1$ and $G_2$. 
Each of $x'_3$ and $x'_4$ sends at most $|G_1|$ edges of type $2$ towards $G_1$, respectively. All other edges between $G_1$ and $G_2$ are the edges from $G_2$ to $G_1$ that are of type $3$. Hence, they are from vertices $u$ where $u^*$ is the leftmost vertex of $G_2$ in its level. 
As the leftmost vertex in each level has at most four children, those vertices in level $\ell'$ send at most $4 a_{\ell'-1,2}+4$ edges towards $G_1$. 
Thus, total number of edges within $D_x$ is at most
\begin{eqnarray*}
e(D_x)&\leq&  e(G_1)+e(G_2)+ e(V(G_1),V(G_2)) \\
&\leq& \sum_{i\in [2]} \left(\frac{14}{3}|G_i|\log_3|G_i| +8 |G_i|\right) +2|G_1|+ \sum_{\ell< \ell'\leq k} (4a_{\ell'-1,2}+4) \\
& \leq& \left(\frac{14}{3}\left(|D_x|+1\right) \log_3|D_x| + 8\left(|D_x|+1\right)\right) + 2|D_x| + 12|D_x|\\
& \leq &\frac{14}{3}|D_x| \log_3|D_x| + 40|D_x|.
\end{eqnarray*}
The penultimate inequality holds since $|D_x|=|G_1|+|G_2|-1$ holds and  we know $a_{\ell'-1,2} = 2(3^{k-\ell'+1}-1)$ and $ \sum_{\ell< \ell'\leq k} 8\cdot 3^{k-\ell'+1} \leq  12\cdot 3^{k-\ell} \leq 12|D_x|$.
\end{description}

Now we count the edges from $D_x$ to $U\setminus D_x$. The edges from $D_x$ to $U\setminus D_x$ must be from a vertex $w \in D_x$ towards $D_T[l(w^*)]$ where $l(w^*)\notin D_x$. Similarly as before, this happens only when $w=x$ or $w^*$ is the leftmost vertex in $D_x$ on a level $\ell' \geq \ell$. 
If $w=x$, it sends at most $a_{\ell-1,2}+1 \leq 6|D_x|$ edges towards $U\setminus D_x$. If $w\neq x$, then there are at most four children of leftmost vertex $w^*$ in a level $\ell'\geq \ell$, the number of such edges is at most 
\[ 6|D_x| + \sum_{\ell'\geq \ell} 4 (a_{\ell', 2}+1) \leq 6|D_x|+ 12\cdot 3^{k-\ell} \leq  18 |D_x|.\]

Hence the total number of edges in $\overrightarrow{G}$ is 
\begin{eqnarray*}
e(\overrightarrow{G}) &\leq& e(U\setminus D[X], U) + \sum_{x\in X}\left(e(D_x,D_x) + e(D_x,U\setminus D_x)\right) \\
&\stackrel{\eqref{eq: edges in U1}}{\leq }& 28 n + \sum_{x\in X} \left(\frac{14}{3}|D_x| \log_3|D_x| + 40 |D_x|\right) + \sum_{x\in X} 18 |D_x| \\
&\leq & \frac{14}{3} n \log_3 n + 100 n.
\end{eqnarray*}
Here, the final inequality holds since $|\sum_{x\in X}D_x|\leq n$.

By \Cref{strong}, $G^2_{T_k}[U]$ is an $n$-vertex universal graph and $G^2_{T_k}[U]$ contains at most $2\cdot 4^2 n \leq 100 n$ more edges than $G_{T_{k}}[U]$ from \Cref{newdef}~\ref{G4}. This is because every vertex in $T_k$ has at most four children. Finally, since $G_{T_k}$ has at most as many edges as $\overrightarrow{G}_{T_k}$, we conclude that $G^2_{T_{k}}[U]$ is a desired universal graph with at most $\frac{14}{3} n \log_3 n + 200 n$ edges. This proves \Cref{thm2}.

\section{Interval-universality}\label{sec:6}

Lemma~\ref{strong} yields an embedding $\phi$ of an arbitrary tree $T'$ with $|T'|\leq |T|$ into our graph $G^{r}_T$. Moreover, it ensures that $V(G^{r}_T)\setminus\phi(T')$ is a $T$-admissible set. In other words, our graph $G^r_T$ is not only universal, but also its induced subgraph $G^r_T[U]$ is universal for any interval $U=\{x_i,\dots, x_j\}$ in the DFS preorder traversal $x_1,\dots, x_n$. This motivates the following new concept.

\begin{definition}
A graph $G$ is an {\bf interval-universal graph} for trees if there exists an ordering $x_1,\dots, x_n$ of the vertices of $G$ in such a way that for any $i,m\in \mathbb{N}$, the graph $G[\{x_{i+1},\dots, x_{i+m}\}]$ contains all $m$-vertex trees as a subgraphs.
\end{definition}

Define $s^{int}(n)$ as the minimum number of edges in an $n$-vertex interval-universal graph for trees. 
\Cref{strong} proves that $G^{2}_T$ for $(2,1)$-tree $T$ is interval-universal with respect to the DFS preorder traversal. Hence we obtain the following inequalities.
\begin{align*}
n\ln{n}-O(n)\leq s^{*}(n)\leq s(n)\leq s^{int}(n)\leq\frac{14}{3}n\log_3{n}+O(n).
\end{align*}
Although improving the lower bounds for $s^{*}(n)$ and $s(n)$ seems very difficult, we can find a larger lower bound for $s^{int}(n)$.

\begin{theorem}\label{thm3}
$s^{int}(n)\geq\left(1-o(1)\right)n\log_2{n}$.
\end{theorem}
\begin{proof}
Let $\varepsilon>0$ be a sufficiently small number.
We use induction on $n$ to prove that $s^{int}(n)\geq (1-\varepsilon) n\log_2 n - \frac{n}{\varepsilon}$.
This is trivially true if $n<1/\varepsilon$, since $\varepsilon$ is small enough.
Assume $n\geq 1/\varepsilon$ and $s^{int}(n') \geq (1-\varepsilon) n'\log_2 n' - \frac{n'}{\varepsilon}$ for all $n'<n$.

It is easy to see that an $n$-vertex universal graph for spanning trees must contain a vertex of degree $n-1$. 
If $G$ is an interval-universal graph for trees and $x_1,\dots, x_n$ is our ordering, then it contains a vertex $x_t$ of degree $n-1$.
Then $G[ \{x_1,\dots, x_{t-1}\}]$ and $G[\{x_{t+1},\dots, x_n\}]$ are both interval-universal, thus induction hypothesis yields 
\begin{align*}
    e(G)&\geq (1-\varepsilon)(t-1)\log_2(t-1)- \frac{t-1}{\varepsilon} + (1-\varepsilon) (n-t)\log(n-t) - \frac{n-t}{\varepsilon} +(n-1) 
    \\
    &\geq  (1-\varepsilon) (n-1)\log_2(\frac{n-1}{2}) - \frac{n-1}{\varepsilon} +(n-1) \geq (1-\varepsilon)n\log_2 n - \frac{n}{\varepsilon}.
\end{align*}
The second inequality is obtained from Jensen's inequality, and the final inequality follows from the AM-GM inequality $\varepsilon (n-1) + \frac{1}{\varepsilon} \geq 2\sqrt{n-1} > \log_2(e)+ \log_2(n)$ where $\varepsilon$ is small and  $n>1/\varepsilon$.
\end{proof}

Note that the gap between the above lower bound and the upper bound from \Cref{thm2} is at most up to the multiplicative constant $\frac{14}{3}\log_32 \leq 2.945$. A natural problem is to improve the upper bound.
One natural idea is to increase the imbalance $K$ to a value larger than two as explained in \Cref{sec: overview} without adding edges from $u$ to the descendants of its nearest-right cousin. Although this will reduce the number of edges, proving that such a graph is universal will need new ideas.

Note that the strategy of generating a graph $G^{r}_T$ from a given tree $T$ and using admissibility in an induction argument yields a bound on $s^{int}(n)$. Thus, just as \Cref{thm3} distinguishes $s^{int}(n)$ from $s^*(n)$ or $s(n)$, some new ideas that only apply to $s^*(n)$ or $s(n)$, but not necessarily to $s^{int}(n)$, might be helpful to further improve upper bounds on $s^*(n)$ or $s(n)$. 

\begin{problem}
Is there a constant $c<\frac{14}{3\ln{3}}$ that satisfies $s^*(n)\leq (c+o(1))n\ln{n}$? 
\end{problem}

\bibliography{references}

@online{chung,
  author       = {Chung, F.},
  title        = {Separator theorems and their applications.},
  url          = {https://mathweb.ucsd.edu/~fan/mypaps/fanpap/117separatorthms.pdf},
}

@article{kaul,
  title={On universal graphs for trees and treewidth $k$ graphs},
  author={Kaul, N. and Wood, D.},
  journal={arXiv preprint arXiv:2508.0335},
  year={2025}
}

@article{alon,
title = {Sparse universal graphs},
journal = {Journal of Computational and Applied Mathematics},
volume = {142},
number = {1},
pages = {1-11},
year = {2002},
author = {Alon, N and Asodi, V.},

}

@article{alonsecond,
author = {Alon, N. and Capalbo, M.},
title = {Sparse universal graphs for bounded-degree graphs},
journal = {Random Structures \& Algorithms},
volume = {31},
number = {2},
pages = {123-133},
year = {2007}
}

@article{szemeredi,
author = {Alon, N. and Capalbo, M. and Kohayakawa, Y. and Rödl, V. and Ruciński, A. and Szemerédi, E.},
year = {2001},
pages = {170-180},
title = {Near-optimum Universal Graphs for Graphs with Bounded Degrees},
volume = {5},
journal = {Approximation, Randomization, and Combinatorial Optimization: Algorithms and Techniques},
}

@article{alonthird,
author = {Alon, N. and Krivelevich, M. and Sudakov, B.},
year = {2007},
pages = {170-180},
title = {Embedding nearly-spanning bounded degree trees},
volume = {27},
number = {6},
journal = {Combinatorica}
}

@article{babai,
  title={On graphs which contain all sparse graphs},
  author={Babai, L. and Chung, F. and Erd{\"o}s, P. and Graham, R. and Spencer, J.},
  journal={North-Holland Mathematics Studies},
  volume={60},
  pages={21-26},
  year={1982}
}

@article{bhatt,
author = {Bhatt, S. and Chung, F. and Leighton, F and Rosenberg, A.},
title = {Universal Graphs for Bounded-Degree Trees and Planar Graphs},
journal = {SIAM Journal on Discrete Mathematics},
volume = {2},
number = {2},
pages = {145-155},
year = {1989},
}

@article{bucic,
title={Universal and unavoidable graphs},
volume={30}, 
number={6}, 
journal={Combinatorics, Probability and Computing}, 
author={Bucić, M. and Draganić, N. and Sudakov, B.}, 
year={2021}, 
pages={942–955}
}

@article{capalbo,
title={Small Universal Graphs for Bounded-Degree Planar Graphs},
volume={22},
number={3},
journal={Combinatorica},
author={Capalbo, M.},
year={2002},
pages={345-359}
}

@article{capalbosecond,
author = {Capalbo, M. and Kosaraju, S.},
title = {Small universal graphs},
year = {1999},
journal = {Proceedings of the Thirty-First Annual ACM Symposium on Theory of Computing},
pages = {741–749}
}

@article{grahamfirst,
author = {Chung, F. and Graham, R.},
title = {On Graphs Which Contain All Small Trees},
year ={1978}, 
pages = {14-23},
volume = {24}, 
journal = {Journal of Combinatorial Theory, Series B}
}

@article{grahamsecond,
author = {Chung, F. and Graham, R.},
title = {ON UNIVERSAL GRAPHS},
journal = {Annals of the New York Academy of Sciences},
year = {2006}, 
volume = {319},
number = {1},
pages = {136-140},
}

@article{grahamthird,
author = {Chung, F. and Graham, R.},
title = {On Universal Graphs for Spanning Trees},
journal = {Journal of the London Mathematical Society},
volume = {s2-27},
number = {2},
pages = {203-211},
year = {1983}
}

@article{pipenger,
author = {Chung, F. and Graham, R. and Pippenger, N.},
title = {On graphs which contain all small trees, \uppercase\expandafter{\romannumeral2}},
year ={1978}, 
pages = {213-223},
volume = {18},
journal = {Colloquia Mathematica Societatis János Bolyai}
}

@article{joret,
author = {Esperet, L. and Joret, G. and Morin, P.},
title = {Sparse universal graphs for planarity},
journal = {Journal of the London Mathematical Society},
volume = {108},
number = {4},
pages = {1333-1357},
year = {2023}
}

@article{gyori,
author = {Gy\H{o}ri, E. and Li, B. and Salia, N. and Tompkins, C.},
title = {A note on universal graphs for spanning trees},
year = {2025},
volume = {362},
journal = {Discrete Applied Mathematics},
pages = {146–147},
}

@article{montgomery,
author = {Montgomery, R.},
title = {Spanning trees in random graphs},
year ={2019}, 
volume = {356},
journal = {Advances in Mathematics},
}
\bibliographystyle{plainurl}

\end{document}